\documentclass[11pt]{article}

\usepackage{times}
\usepackage{amsfonts,amssymb}
\usepackage[a4paper]{geometry}
\newtheorem{thm}{Theorem}[section]
\newtheorem{lem}[thm]{Lemma}
\newenvironment{proof}{\noindent {\bf Proof}.}{\hfill$\Box$\\[-5mm]}

\newcommand{\NN}{\mathbb{N}}
\newcommand{\RR}{\mathbb{R}}
\newcommand{\DD}{\mathbb{D}}

\newcommand{\EE}{\mathbf{E}}
\newcommand{\PP}{\mathbf{P}}

\newcommand{\One}{{\mathrm{1} \kern -0.27em \mathrm{I}}}

\newcommand{\defeg}{\mathrel :=}

\newcommand{\ut}{{\tilde u}}
\newcommand{\lin}{\mathbf{K}}
\newcommand{\lint}{\mathcal{K}}

\begin{document}

\title{Functional central limit theorems 
for a large network in which customers join the shortest of several queues}

\author{\sc Carl Graham~\footnote{CMAP, 
{\'E}cole Polytechnique, 91128 Palaiseau, France. UMR CNRS 7641. 
\tt carl@cmapx.polytechnique.fr}
}
\date{}

\maketitle

\begin{small}
\noindent
{\bf Abstract}.
\addtolength{\baselineskip}{.2\baselineskip}
We consider
$N$ single server infinite buffer queues with service rate $\beta$.
Customers arrive at rate $N\alpha$,
choose $L$ queues uniformly, and join the shortest.
We study the processes 
$t\in\RR_+ \mapsto R^N_t=(R^N_t(k))_{k \in\NN}$ for large $N$, where
$R^N_t(k)$ is the fraction of queues of length at least $k$ at time $t$.
Laws of large numbers (LLNs) are known, see
Vvedenskaya et al.~\cite{Vved:96}, Mitzenmacher~\cite{Mitzenmacher:96}
and Graham \cite{Graham:00a}.
We consider certain Hilbert spaces with the weak topology.
First, we prove a
functional central limit theorem (CLT)
under the {\em a~priori\/} assumption that the initial data 
$R^N_0$
satisfy the corresponding CLT. We use a compactness-uniqueness method,
and the  limit is characterized as an Ornstein-Uhlenbeck (OU) process.
Then, we study 
the $R^N$ in equilibrium under the stability  condition $\alpha < \beta$,
and prove a functional CLT 
with limit the OU~process in equilibrium.
We use ergodicity  and
justify the  inversion of limits 
$\lim_{N\to\infty}\lim_{t \to\infty} =\lim_{t\to\infty}\lim_{N \to\infty}$
by a compactness-uniqueness method.
We deduce {\em a posteriori\/} 
the CLT for $R^N_0$ under the invariant laws, an interesting result in its own right.
The main tool for proving tightness of the implicitly defined 
invariant laws in the CLT scaling
and ergodicity of the limit OU~process
is a global  exponential stability result
for the nonlinear dynamical system obtained in the functional LLN
limit.



\noindent
March 25, 2004.
\end{small}
\bigskip\hrule

\addtolength{\baselineskip}{.4\baselineskip}

\section{Introduction}

\subsection{Preliminaries}

We consider a Markovian  network constituted of $N\ge L\ge1$ 
infinite buffer single server queues.
Customers arrive at rate $N\alpha$, are each allocated
$L$ distinct queues uniformly at random, and join the shortest, ties
being resolved uniformly. 
Servers work at rate $\beta$.
Arrivals, allocations, and services 
are independent.
For $L=1$ we have i.i.d.\ $M_\alpha/M_\beta/1/\infty$ queues. 
For $L\ge2$ the interaction structure depends only on sampling 
from the empirical measure of $L$-tuples of queue
states:
in statistical mechanics terminology, the system is in 
$L$-body mean-field interaction.
We continue the large $N$ study  introduced
by Vvedenskaya et al.~\cite{Vved:96} 
and Mitzenmacher~\cite{Mitzenmacher:96} and 
continued in Graham~\cite{Graham:00a}.

The process $(X^N_i)_{1\leq i\leq N }$ is Markov, where
$X^N_i(t)$ denotes the length of queue $i$ at time $t$ in $\RR_+$. 
Its empirical measure $\mu^N  ={1 \over N}\sum_{i=1}^N \delta_{X^N_i}$ 
has samples in $\mathcal{P}(\DD(\RR_+,\NN) )$, and its marginal process 
$\bar X^N = (\bar X^N_t)_{t\ge0}$ with
$\bar X^N_t = \mu_t={1 \over N}\sum_{i=1}^N \delta_{X^N_i(t)}$
has sample paths  
in $\DD\!\left(\RR_+,\mathcal{P}(\NN)\right)$. 
We are interested 
in the tails of the marginals $\bar X^N_t$ and 
consider
\[
R^N =(R^N_t)_{t\ge0}\,,
\qquad
R^N_t = (R^N_t(k))_{k\in\NN}\,,
\qquad
R^N_t(k) = {1 \over N}\sum_{i=1}^N \One_{X^N_i(t)\ge k}\,,
\] 
and $R^N_t(k)$ is the fraction of queues of length at least $k$
at time $t$.
For the uniform topology on
\[
\mathcal{V} 
= \left\{ 
(v(k))_{k\in\NN} : v(0)=1,\ v(k)\geq v(k+1),\ \lim v =0 
\right\} \subset c_0\,,
\qquad 
\mathcal{V}^N = \mathcal{V} \cap {1\over N}  \NN^\NN\,,
\]
coinciding here with the product topology,
the process $R^N$ has sample paths in $\DD\!\left(\RR_+,\mathcal{V}^N\right)$.

The processes $\bar X^N$ and $R^N$ 
are in relation through
$p \in \mathcal{P}(\NN) \longleftrightarrow v \in \mathcal{V}$
for $v(k) = p[k,\infty )$ and $p\{k\} = v(k) - v(k+1)$ for
$k$ in $\NN$. This classical homeomorphism 
maps the subspace of probability measures with 
finite first moment
onto  $\mathcal{V}\cap \ell_1$,  corresponding to
a finite number of customers in the network.
The symmetry structure implies that these processes 
are Markov.

The stationary regime has great practical relevance.
The stability condition $\alpha < \beta$
(Theorem~5~(a) in \cite{Vved:96}, Lemma 3.1 in \cite{Mitzenmacher:96},
Theorem~4.2 in \cite{Graham:00a}) is obtained from ergodicity criteria 
yielding little information.
We study the large $N$ asymptotics of $R^N$, first for
transient regimes with appropriately converging initial data, and then 
in equilibrium using an
indirect approach involving ergodicity in well-chosen transient regimes
and an inversion of limits 
for large $N$ and large times.
Law of large numbers (LLN) results are already known, and
we obtain functional central limit theorems (CLTs).

\subsection{Previous results: laws of large numbers}

We relate results found in essence in Vvedenskaya et al.~\cite{Vved:96}. 
We follow Graham~\cite{Graham:00a} which extends these results, notably by
considering  the  empirical measures on path space 
$\mu^N$ and thus yielding chaoticity results (asymptotic independence of queues).
Chapter 3 in Mitzenmacher~\cite{Mitzenmacher:96} gives related results.
(The rates $\alpha$ and $\beta$ correspond to
$\lambda$ and $1$ in  \cite{Vved:96,Mitzenmacher:96} and 
$\nu$ and $\lambda$ in \cite{Graham:00a}.)

Consider the  mappings with values in $c^0_0$ given
for $v$ in $c_0$ by
\begin{equation}
\label{F}
F_+(v)(k) = \alpha\!\left( v(k-1)^L -v(k)^L  \right),
\quad 
F_-(v) (k) = \beta  (v(k) - v(k+1))\,,
\qquad 
k\ge1\,,
\end{equation}
and $F = F_+ - F_-$, and the nonlinear differential equation 
$\dot u = F(u)$
on $\mathcal{V}$, given for $t\ge0$ by 
\begin{eqnarray}
\label{is}
\dot u_t(k)
&=& \alpha\! \left( u_t(k-1)^L-u_t(k)^L \right)
- \beta\! \left( u_t(k) - u_t(k+1)\right)
\nonumber\\
&=& 
\alpha u_t(k-1)^L
- \left(\alpha u_t(k)^L + \beta u_t(k) \right) 
- \beta u_t(k+1)\,,
\qquad k\geq 1\,.
\end{eqnarray}
This corresponds to the systems
(1.6) in \cite{Vved:96}, (3.5) in  \cite{Mitzenmacher:96}  
and  (3.9) in  \cite{Graham:00a}. Note that $F_-$ is linear.

\begin{thm}
\label{euv}
There  exists a unique solution $u=(u_t)_{t\ge0}$
taking values in $\mathcal{V}$ for (\ref{is}), 
and $u$ is in $C(\RR_+,\mathcal{V})$. 
If $u_0$ is in $\mathcal{V}\cap\ell_1$
then $u$ takes values in $\mathcal{V}\cap\ell_1$.
\end{thm}

\begin{proof}
We use Theorem~3.3  and Proposition 2.3 in \cite{Graham:00a}. 
These exploit the homeomorphism
between $\mathcal{P}(\NN)$ with the weak topology 
and  $\mathcal{V}$ with the product topology.
Then (\ref{is}) corresponds to a non-linear forward 
Kolmogorov equation for a pure jump
process with uniformly bounded (time-dependent) jump rates.
Uniqueness within the class of bounded measures
and existence of a probability-measure valued solution 
are obtained using the total variation norm.
Theorem~1~(a) in \cite{Vved:96} yields existence (and uniqueness) 
in $\mathcal{V}\cap\ell_1$.
\end{proof}

Firstly, a functional LLN
for initial conditions satisfying a LLN
is part of Theorem~3.4 in \cite{Graham:00a} and can be deduced from
Theorem~2 in \cite{Vved:96}.

\begin{thm}
\label{lln}
Assume that
$(R^N_0)_{N\ge L}$ converges in law  to
$u_0$ in $\mathcal{V}$.
Then  $(R^N)_{N\ge L}$ converges in law in $\DD(\RR_+,\mathcal{V})$
to the unique solution $u=(u_t)_{t\ge0}$ 
starting at $u_0$ for (\ref{is}).
\end{thm}

Secondly, for $\alpha<\beta$ 
the limit equation  (\ref{is}) has a globally attractive stable point $\ut$ 
in $\mathcal{V}\cap \ell_1$. 

\begin{thm}
\label{utilde}
Let $\rho = \alpha / \beta <1$. The  equation (\ref{is}) has
a unique stable point in $\mathcal{V}$ given by 
\[
\ut = (\ut(k)_{k\in\NN}\,,
\qquad
\ut(k)= \rho^{(L^k -1) / (L-1)}= \rho^{L^{k-1}+L^{k-2}+\cdots + 1}\,,
\]
and the solution $u$ of (\ref{is}) 
starting at any $u_0$ in $\mathcal{V}\cap \ell_1$  is such that
$\lim_{t\to\infty} u_t = \ut$.
\end{thm}
\begin{proof}
Theorem~1~(b) in \cite{Vved:96} yields that $\ut$ 
is  globally asymptotically stable in $\mathcal{V}\cap \ell_1$.
A stable point $u$ in $\mathcal{V}$ satisfies
$
\beta u(k+1) - \alpha u(k)^L = \beta u(k) - \alpha u(k-1)^L
= \cdots =
\beta u(1) - \alpha 
$
and converges to $0$, hence
$u(1) = \alpha/\beta$ and
$u(2)$,  $u(3)$, \dots\, are successively determined uniquely.
\end{proof}

Lastly, 
a compactness-uniqueness argument
justifies the inversion of limits 
$\lim_{N\to\infty}\lim_{t\to\infty} = \lim_{t\to\infty}\lim_{N\to\infty}$,
which yields a result in equilibrium. This 
method,
used by Whitt~\cite{Whitt:85} for the star-shaped loss network,
is detailed in Graham~\cite{Graham:00b} Sections~9.5 and 9.7.3.
The following functional LLN in equilibrium (Theorem~4.4 in \cite{Graham:00a})
can be deduced from \cite{Vved:96} but is not stated there, and
implies that under the invariant laws
$\lim_{N\to\infty} \EE( R^N_0(k) ) = \ut(k)$ for $k\in\NN$
(Theorem~5 (c) in \cite{Vved:96}).

\begin{thm}
\label{lln.eq}
Let $\rho = \alpha / \beta <1$ and
the networks of size $N$ be in equilibrium. Then $(R^N)_{N\ge L}$
converges in probability in $\DD(\RR_+,\mathcal{V})$ to $\ut$.
\end{thm}

Note that $\ut(k)$ decays hyper-exponentially in $k$ for $L\ge2$ 
instead of the exponential decay $\rho^{k}$ corresponding to 
i.i.d.\ queues in equilibrium ($L=1$). 
For finite networks in equilibrium there is at most exponential decay since
$
\PP\!\left(X^N_1 + \cdots + X^N_N \ge N k\right) \le 
\PP\!\left(X^N_1\ge k\right) + \cdots + \PP\!\left(X^N_N \ge k\right)
$
and 
by comparison with an $M_{N\alpha}/M_{N\beta}/1$ queue
\begin{equation}
\label{finexp}
\EE\!\left(R^N_t(k)\right) 
= \PP\!\left(X^N_i(t) \ge k\right) \ge {1\over N} \rho^{Nk}\,,
\qquad
k\ge0\,.
\end{equation}

The asymptotic queue sizes are dramatically decreased 
by this simple load balancing (or resource pooling) procedure,
which carries little overhead even for large $N$ since $L$
is fixed (for instance $L=2$). 
This feature is quite robust and
true for many systems, as
was illustrated on several examples by 
Mitzenmacher~\cite{Mitzenmacher:96} and Turner~\cite{Turner:98} 
using proofs as well 
as simulations.
It can be used as a guideline for
designing practical networks. 
In contrast, the bound (\ref{finexp}) assumes
the best utilization of the $N$ servers, fully collaborating even for a single customer.

Theorem~3.5 in  Graham~\cite{Graham:00a} gives convergence bounds
on bounded time intervals $[0,T]$ for  i.i.d.  $(X^N_i(0))_{1\leq i \leq N}$
using results in Graham and M{\'e}l{\'e}ard~\cite{Meleard:94}. This can be extended 
if the initial laws satisfy {\em a~priori} controls, 
but it is not so in equilibrium (the bounds are exponentially large in $T$).

\subsection{The outline of this paper}

The study of the fluctuations around the functional LLN will yield  for instance
asymptotically tight confidence intervals  for the process
$t \mapsto N^{-1}\mathrm{Card}\!\left\{i=1,\ldots,N : X^N_i(t) \in A \right\}$.
In a realistic setting (finite number of finite buffer queues) 
such confidence intervals would allow network evaluation 
or dimensioning in function of quality of 
service requirements on delays and overflows.  
The LLN on path space concerns objects such as
$N^{-1}\mathrm{Card}\!\left\{i=1,\ldots,N 
: \left(t \mapsto X^N_i(t)\right) \in B \right\}$
with a richer temporal structure, 
but topological difficulties usually block the corresponding fluctuation study.

We consider the process $R^N$ and solution $u$ for (\ref{is})
starting at $R^N_0$ in  $\mathcal{V}^N$ and $u_0$ in $\mathcal{V}$, and
\begin{equation}
\label{flu}
 Z^N = \sqrt{N} (R^N - u).
\end{equation}
The processes $Z^N = (Z^N_t)_{t\ge0}$
will be studied in the Skorokhod spaces on appropriate
Hilbert spaces with the weak topology. These spaces are not metrizable
and require appropriate tightness criteria.

We first consider a wide class of 
$R^N_0$ and $u_0$ under the {\em assumption\/} that
$(Z^N_0)_{N\ge L}$ converges in law (for instance satisfies a CLT). 
We obtain
a functional CLT  
in relation to Theorem~\ref{lln}, with limit given by an
Ornstein-Uhlenbeck (OU) process starting at the limit of the $(Z^N_0)_{N\ge L}$.
This covers without constraints on $\alpha$ and $\beta$ many {\em transient\/} regimes
with {\em explicit\/} initial conditions, such as 
initially empty networks, or more generally
i.i.d.\ initial queue sizes.

We then focus on  the {\em stationary\/} regime for $\alpha<\beta$.
The initial data is now {\em implicit\/}:
the law of $R^N_0$ is the invariant law for $R^N$
and $u_0=\ut$.
We prove tightness for  $(Z^N_0)_{N\ge L}$ 
using the ergodicity of $Z^N$ for fixed $N$
and intricate fine studies of the long-time behavior of the  
nonlinear dynamics appearing at the 
large $N$ limit.
The main result in this paper is
a functional CLT
in equilibrium for
$(Z^N)_{N\ge L}$ with limit the OU~process in equilibrium.
This {\em implies\/} a CLT under the
invariant laws for $(Z^N_0)_{N\ge L}$,  
an important result which seems difficult to obtain directly.

Section~\ref{clt} introduces without proof the main notions and results. 
Section~\ref{Sproofs} gives the proof of the functional CLT
for converging initial data
by compactness-uniqueness and martingale techniques.

We then consider $u_0=\tilde u$.
We study the OU~process in Section~\ref{propOU},
derive a spectral representation for the linear operator in the drift,
and prove the existence of a spectral gap. 
A main difficulty is that the scalar product for which the operator is self-adjoint
is {\em too strong\/}
for the limit dynamical system and the invariant laws for finite $N$.
We consider appropriate Hilbert spaces in which the operator 
is {\em not\/} self-adjoint
and prove exponential stability.

In Section~\ref{expstab} we likewise prove that $\ut$ is globally
exponentially stable for the non-linear 
dynamical system.
In Section~\ref{cfp} we obtain
bounds for  $Z^N_t$ uniform for $t\ge0$ and large  $N$,
using  the preceding stability result in order 
to  iterate the bounds on intervals of length $T$. 
Bounds on the invariant laws of $Z^N$ follow using ergodicity.
The proof for the functional CLT in equilibrium follows from
a compactness-uniqueness argument involving the  functional CLT for 
converging initial data.

\section{The functional central limit theorems}
\setcounter{equation}{0}
\label{clt}

\subsection{Preliminaries}

The exponential of a bounded linear operator is given by the usual series expansion.
Let $c^0_0$ and $\ell^0_p$ for $p\ge1$ be
the subspaces  of sequences vanishing at $0$
of the classical sequence spaces $c_0$ (with limit 0) and
$\ell_p$ (with summable $p$-th power).
In matrix notation we use
the canonical basis, hence
sequences vanishing at $0$
are identified with infinite column vectors
indexed by $\{1,2,\cdots\}$.
The diagonal matrix with terms given by
the sequence $a$ is denoted by $\mathrm{diag}(a)$. 
Sequence inequalities, etc., should be interpreted termwise.
Empty sums are equal to
$0$ and empty products to $1$. Constants such as $K$ may
vary from line to line.
Let $g_\theta = (\theta^k)_{k\ge1}$ be the geometric sequence with parameter $\theta$.

For a sequence $w=(w(k))_{k\ge 1}$ such that $w(k)>0$
we define the Hilbert spaces
\[
L_2(w) = \Biggl\{ x \in \RR^\NN : x(0)=0\,,\ 
\Vert x \Vert_{L_2(w)}^2 = \sum_{k\ge1} \biggl({x(k)\over w(k)}\biggr)^2 w(k)
= \sum_{k\ge1} x(k)^2 w(k)^{-1} <\infty
\Biggr\}
\]
and in matrix notation $(x,y)_{L_2(w)} = x^* \mathrm{diag}(w^{-1})y$.
We use the notation $L_2(w)$ since its 
elements will often be considered as measures identified with 
their densities with respect to the reference measure $w$. 
In this perspective $L_1(w) = \ell_1^0$ and if
$w$ is summable then
$\Vert x \Vert_1 \le \Vert w \Vert_1^{1/2}\Vert x \Vert_{L_2(w)}$
and $L_2(w) \subset \ell_1^0$. Using $L_2(1) = \ell^0_2$ as a pivot space,
for bounded $w$ we have the Gelfand triplet
$
L_2(w) \subset \ell_2^0  \subset L_2(w)^* = L_2(w^{-1})
$.

Another natural perspective on $L_2(w)$ is that it is
an $\ell_2$ space with weights, and 
we consider the $\ell_1$ space with same weights
(the notation being chosen for consistency)
\[
\ell_1(w) = \Biggl\{ x \in \RR^\NN : x(0)=0\,,\ 
\Vert x \Vert_{\ell_1(w)} = \sum_{k\ge1} |x(k)| w(k)^{-1} <\infty
\Biggr\}
\]
and $x \in L_2(w) \Leftrightarrow x^2 \in \ell_1(w)$ with  
$\Vert x \Vert_{L_2(w)}^2 = \Vert x^2 \Vert_{\ell_1(w)}$.
The inclusion $\mathcal{V} \cap \ell_1(w) \hookrightarrow \mathcal{V} \cap L_2(w)$
is continuous since $x^2 \le |x|$ for $|x| \le 1$.
The following result is trivial.

\begin{lem}
\label{eqiv.l2w}
If $w = O(v)$ and $v = O(w)$
then the $L_2(v)$ and  $L_2(w)$ norms
are equivalent, and the  $\ell_1(v)$ and  $\ell_1(w)$ norms
are equivalent.
\end{lem}

In the sequel we often assume that $w=(w_k)_{k\ge1}$ satisfies the condition that 
\begin{equation}
\label{compw}
\exists\, c,d>0\,,\, \forall k\ge1\,,\, 
0< c  w(k+1)\le  w(k) \le d w(k+1)\,,
\end{equation}
which is satisfied by $g_\theta=(\theta^k)_{k\ge1}$ with
$c=d=1/\theta$ for $\theta > 0$.
It implies that $ w(1) d (1/d)^k \le w(k) \le w(1) c (1/c)^k$
which bounds $w$ by geometric sequences. 
The norms have exponentially strong weights for $c>1$.
We give a refined existence result for $(\ref{is})$.  
(Proofs are left for later.)

\begin{thm}
\label{isw}
Let $w$ satisfy (\ref{compw}).
Then in $\mathcal{V}$ the mappings
$F$, $F_+$ and $F_-$ are Lipschitz for 
the $L_2(w)$ and the $\ell_1(w)$ norms. Existence and uniqueness holds
for (\ref{is}) in  $\mathcal{V} \cap L_2(w)$ and in $\mathcal{V} \cap \ell_1(w)$.
\end{thm}

\subsection{The functional CLT for converging initial data}

\subsubsection*{The time-inhomogeneous Ornstein-Uhlenbeck process}

In $\mathcal{V}$, the linearization of (\ref{is}) around a particular solution $u$
is the linearization of the recentered equation
satisfied by $y=g -u$ 
where $g$ is a generic solution for (\ref{is}). 
It is given for $t\ge0$ by
\begin{equation}
\label{linis.gen}
\dot z_t = \lin(u_t) z_t
\end{equation} 
where for $v$ in $\mathcal{V}$ the linear operator 
$\lin(v) : x \mapsto \lin(v) x$ on $c_0^0$ is given by
\begin{equation}
\label{defK}
\lin(v) x (k) =
\alpha L v(k-1)^{L-1} x(k-1)
-(\alpha L v(k)^{L-1} + \beta) x(k)
+  \beta x(k+1)\,,\quad 
k\ge1\,,
\end{equation}
and is identified with its infinite matrix in the canonical basis 
$(0,1,0,0 \ldots)$, $(0,0,1,0 \ldots)$, \dots\, 
\\
\[
\lin(v)=
\pmatrix{
 -\left(\alpha L v(1)^{L-1} + \beta\right) &  \beta & 0  &\cdots 
\cr
 \alpha L v(1)^{L-1} & -\left( \alpha L v(2)^{L-1} + \beta\right) &  \beta   & \cdots
\cr
 0 &  \alpha L v(2)^{L-1} & -\left( \alpha L v(3)^{L-1} + \beta\right)  & \cdots
\cr
 0 &  0 & \alpha L v(3)^{L-1} & \cdots
\cr
 \vdots & \vdots & \vdots  & \ddots 
}.
\]

Let  $(M(k))_{k\in\NN}$ be independent
real continuous centered Gaussian martingales, determined in law by their
deterministic Doob-Meyer brackets given for $t\ge 0$ by
\begin{equation}
\label{limdoob}
 \langle M(k) \rangle_t = \int_0^t
\left\{ F^{\vphantom{N}}_+(u_s)(k) + F^{\vphantom{N}}_-(u_s)(k) \right\}ds\,. 
\end{equation}
The processes $M = (M(k))_{k\ge0}$ and
$\langle M \rangle = \left(\langle M(k) \rangle\right)_{k\in\NN}$
have values in $c^0_0$.

\begin{thm}
\label{gmsi}
Let $w$ satisfy (\ref{compw}) and $u_0$ be in $\mathcal{V} \cap \ell_1(w)$.
Then the Gaussian martingale $M$ is square-integrable in $L_2(w)$.
\end{thm}

\begin{proof}
We have 
$\EE\!\left(\Vert M_t \Vert_{L_2(w)}^2\right)
= \Vert \langle M\rangle_t \Vert_{\ell_1(w)}$
and we conclude using (\ref{limdoob}), Theorem \ref{isw},
and uniform bounds in 
$\ell_1(w)$ on $(u_s)_{0\le s \le t}$  in function of $u_0$
given by the  Gronwall Lemma.
\end{proof}

The limit equation for the fluctuations is 
a Gaussian perturbation of (\ref{linis.gen}), 
the inhomogeneous affine SDE  given for $t\ge0$ by
\begin{equation}
\label{ousde}
Z_t = Z_0 + \int_0^t \lin(u_s) Z_s \,ds + M_t\,.
\end{equation}
A well-defined solution is called an Ornstein-Uhlenbeck
process, in short OU~process.
We recall that strong (or pathwise) uniqueness implies weak uniqueness, and
that $\ell_1(w) \subset L_2(w)$.

\begin{thm}
\label{Kbdd}
Let the sequence  $w$ satisfy (\ref{compw}). 
\\
(a) 
For $v$ in $\mathcal{V}$, 
the operator $\lin(v) $ is bounded in $L_2(w)$ with
 operator norm uniformly bounded in $v$.
\\
(b) Let $u_o$ be in $\mathcal{V} \cap L_2(w)$. Then in $L_2(w)$ there
is a unique solution 
$z_t = \mathrm{e}^{\int_0^t \lin(u_s)\,ds} z_0$ for (\ref{linis.gen}) and
strong uniqueness of solutions holds for (\ref{ousde}).
\\
(c) Let $u_o$ be in $\mathcal{V} \cap\ell_1(w)$.
Then in $L_2(w)$ there
is a unique strong solution 
$
Z_t = \mathrm{e}^{\int_0^t \lin(u_s)\,ds} Z_0 
+ \int_0^t \mathrm{e}^{\int_s^t \lin(u_r)\,dr} dM_s
$
for (\ref{ousde})
and if $\EE\Bigl(\Vert Z_0 \Vert_{L_2(w)}^2\Bigr) < \infty$ then 
$\EE\Bigl(\sup_{t\le T}\Vert Z_t \Vert_{L_2(w)}^2\Bigr) < \infty$.
\end{thm}

\subsubsection*{Tightness bounds and the CLT}

The finite-horizon bounds in the following  lemma  will yield
tightness estimates for the processes $Z^N$
used in the compactness-uniqueness proof for the 
subsequent theorem.

\begin{lem}
\label{l2bds}
Let $w$ satisfy (\ref{compw}).
Let $u_0$ be in $\mathcal{V} \cap\ell_1(w)$
and $R^N_0$ be in $\mathcal{V}^N$.
For any $T\ge0$
\[
\limsup_{N\to\infty} 
\EE\left(\left\Vert Z^N_0 \right\Vert_{L_2(w)}^2\right) < \infty
\Rightarrow
\limsup_{N\to\infty} 
\EE\biggl( \sup_{0 \le t\le T} \left\Vert Z^N_t \right\Vert_{L_2(w)}^2\biggr) 
< \infty\,.
\]
\end{lem}

We refer to Jakubowski~\cite{Jakubowski:86} for the Skorokhod topology
for non-metrizable topologies.
For the weak topology of a 
reflexive Banach space, the relatively compact sets are the 
bounded sets for the norm, see Rudin~\cite{Rudin:73} Theorems 1.15~(b),
3.18, and 4.3.
Hence, if  $B(r)$ denotes the closed ball centered at 0 of radius $r$, 
a set $\mathcal{T}$ of probability measures
is tight if and only if for all $\varepsilon>0$ there exists 
$r_\varepsilon < \infty$ such that $p(B(r_\varepsilon))> 1 - \varepsilon$
uniformly for $p$ in $\mathcal{T}$. We state the functional CLT.

\begin{thm}
\label{fclt}
Let $w$ satisfy (\ref{compw}). Consider 
$L_2(w)$
with its weak topology and
$\DD(\RR_+, L_2(w))$ with the corresponding Skorokhod
topology.
Let $u_0$ be in $\mathcal{V} \cap \ell_1(w)$,
$R^N_0$ in $\mathcal{V}^N$, and
$Z^N$ be given by (\ref{flu}).
If $(Z^N_0)_{N\ge L}$ converges in law to $Z_0$ and is tight,
 then $(Z^N)_{N\ge L}$
converges in law to the unique  OU~process 
solving (\ref{ousde}) starting at $Z_0$ and is tight.
\end{thm}

\subsection{The functional CLT in equilibrium}
\label{oueq}

We assume the stability condition  $\rho = \alpha / \beta <1$ holds,
and consider  $u_0=\tilde u$.

\subsubsection*{The Ornstein-Uhlenbeck process}

We set $\lint = \lin(\tilde u)$ and (\ref{defK}) yields that
$\lint : x \in c_0^0 \mapsto \lint x \in c_0^0$ is given  by
\begin{equation}
\label{defKst}
\lint x(k) = \lin(\tilde u) x(k)
=
\beta L \rho^{L^{k-1}} x(k-1)
-\left(\beta L \rho^{L^k} + \beta\right) x(k)
+  \beta x(k+1)\,,
\quad k\ge1\,,
\end{equation}
identified with its infinite matrix in the canonical basis
\begin{equation}
\label{matrixK}
\lint =
\pmatrix{
 -\left(\beta L \rho^{L} + \beta\right) &  \beta & 0  & 0 &\cdots 
\cr
 \beta L \rho^{L} & -\left( \beta L \rho^{L^2} + \beta\right) &  \beta   & 0 & \cdots
\cr
 0 &  \beta L \rho^{L^2} & -\left( \beta L \rho^{L^3} + \beta\right)  & \beta & \cdots
\cr
 0 &  0 &  \beta L \rho^{L^3} & -\left( \beta L \rho^{L^4} + \beta\right)  &  \cdots
\cr
 \vdots & \vdots & \vdots  & \vdots & 
}.
\end{equation}

Note that $\lint=\mathcal{A}^*$ where $\mathcal{A}$  is  
the generator of a sub-Markovian birth and death process.
We give the Karlin-McGregor spectral decomposition for $\lint$ 
in Section~\ref{spectdec}, to which we make a few forward
references (it is {\em not\/} a resolution of the identity, 
see Rudin~\cite{Rudin:73}). 
The potential coefficients  of $\mathcal{A}$ are 
\begin{equation}
\label{pot}
\pi = (\pi(k))_{k\ge1}\,,
\qquad
\pi(k) =   L^{k-1} \rho^{(L^k-L)/(L-1)} =  \rho^{-1} L^{k-1} {\ut}(k)\,,
\end{equation}
and solve the detailed balance equations $\pi(k+1) = L\rho^{L^k}\pi(k)$
with $\pi(1)=1$. 
The linearization of (\ref{is}) around its stable point $\ut$
is the forward Kolmogorov equation for  $\mathcal{A}$ 
given for $t\ge0$ by
\begin{equation}
\label{linis}
\dot z_t = \lint z_t
\end{equation}
which is special case of (\ref{linis.gen}). 
Considering (\ref{limdoob}) and $F(\tilde u)=F_+(\tilde u) - F_-(\tilde u)=0$, 
the martingale $M=(M(k))_{k\in\NN}$ 
has the same law as a  $c_0^0$-valued sequence
$B=(B(k))_{k\in\NN}$ of independent centered Brownian motions
such that
$B(0)=0$ and for $k\ge1$
\[
\tilde v(k):=
\mathrm{var}(B_1(k)) = \EE\!\left(B_1(k)^2\right) 
= 2 \beta \left(\ut(k) - \ut(k+1)\right)
= 2 \beta  \rho^{(L^k-1)/ (L-1)} \left(1 - \rho^{L^{\smash{k}}} \right),
\]
and $B$ has diagonal
infinitesimal covariance matrix $\mathrm{diag}(\tilde v)$. The following 
result is obvious.

\begin{thm}
\label{chbm}
The process $B$ is an Hilbertian Brownian motion in $L_2(w)$ 
if and only if $\ut$ is in $\ell_1(w)$.
This is true for 
$w= \pi$ and $w = g_\theta$ for $\theta>0$ when $L\ge2$ or 
for $w = g_\theta$ for $\theta>\rho$ when $L=1$.
\end{thm}

The Ornstein-Uhlenbeck (OU) process $Z = (Z(k))_{k\in\NN}$ solves
the affine SDE given for $t\ge0$ by 
\begin{equation}
\label{sde}
Z_t = Z_0 + \int_0^t \lint Z_s \,ds + B_t
\end{equation} 
which is a Brownian perturbation of (\ref{linis}).
For  $L\ge2$, existence and uniqueness results hold under much weaker
assumptions than (\ref{compw}).

\begin{thm}
\label{Dbdd}
 Let $w$ be such that there exists  $c>0$ and $d>0$ with
\[
0 < c  w(k+1)\le  w(k) \le d \rho^{-2L^k} w(k+1)\,,
\qquad k\ge1\,.
\]
(a)
In $L_2(w)$, the operator $\lint$ is bounded,
the equation (\ref{linis}) has a unique solution 
$z_t = \mathrm{e}^{\lint t} z_0$ where
$\mathrm{e}^{\lint t}$  has a spectral representation 
given by (\ref{repr}), and
there is uniqueness of solutions for the SDE~(\ref{sde}).
The assumptions and  conclusions hold 
for $w= \pi$ and $w = g_\theta$ for $\theta>0$.
\\
(b) In addition let $w$ be such that  $\ut$ is in $\ell_1(w)$.
The SDE (\ref{sde}) has a unique solution
$
Z_t = \mathrm{e}^{\lint t} Z_0 
+ \int_0^t \mathrm{e}^{\lint(t-s)}\,dB_s
$ in $L_2(w)$ further made explicit in (\ref{explou}).
This the case  for 
$w= \pi$ and $w = g_\theta$ for $\theta>0$ when $L\ge2$ or
for $w = g_\theta$ for $\theta>\rho$ when $L=1$ .
\end{thm}

We use results in van Doorn~\cite{Doorn:85} to prove the existence
of a spectral gap, and
use this fact for an exponential stability  result inspired from
Callaert and Keilson~\cite{Callaert:73b} Section 10.

\begin{thm}
\label{spe.gap}
(Spectral gap.)
The operator $\lint$ is bounded self-adjoint in $L_2(\pi)$.
The least point $\gamma$ of the spectrum of $\lint$
is such that $0  < \gamma \le \beta$. The 
solution $z_t = \mathrm{e}^{\lint t}z_0$ for (\ref{linis}) in $L_2(\pi)$
satisfies
$\Vert z_t \Vert_{L_2(\pi)} 
\le \mathrm{e}^{-\gamma t}  \Vert z_0 \Vert_{L_2(\pi)}$.
\end{thm}

For $L\ge2$ the sequence $\pi$ decays hyper-exponentially, see (\ref{pot}),
and (\ref{finexp}) implies that the  $L_2(\pi)$ norm is too strong 
for the CLT. 
Further,
the mapping $F_+$ is not Lipschitz in 
$\mathcal{V} \cap L_2(\pi)$ for 
the $L_2(\pi)$ norm,
see Theorems~\ref{isw} and \ref{Dbdd} and their
contrasting assumptions and proofs. Hence we prove exponential
stability and (exponential) ergodicity for the OU~process
in appropriate spaces.

\begin{thm}
\label{linis.exp}
Let  $0 < \theta < 1$ when $L\ge2$ or $\rho \le \theta < 1$ when $L=1$. 
There  exists  $\gamma_\theta>0$ and $C_\theta < \infty$ such that
the solution $z_t = \mathrm{e}^{\lint t}z_0$ for (\ref{linis})  
in $L_2(g_\theta)$
satisfies
$\Vert z_t \Vert_{L_2(g_\theta)} 
\le \mathrm{e}^{-\gamma_\theta t} C_\theta \Vert z_0 \Vert_{L_2(g_\theta)}$.
\end{thm}

\begin{thm}
\label{ou.erg}
Let $w= \pi$ or $w = g_\theta$ with $0 < \theta<1$ when $L\ge2$ or
let $w = g_\theta$ with $\rho < \theta<1$ when $L=1$.
Any solution for the SDE (\ref{sde}) in $L_2(w)$
converges in law for large times 
 to its unique invariant law (exponentially fast). This law is
the law of $\int_0^\infty \mathrm{e}^{\lint t}dB_t$ which is 
 Gaussian centered with covariance matrix 
$\int_0^\infty \mathrm{e}^{\lint t} \mathrm{diag}(\tilde v) \mathrm{e}^{\lint^* t}dt$
made more explicit in  (\ref{limlaw}) and (\ref{cov}).
There is a unique stationary OU~process
solving the SDE (\ref{sde}) in  $L_2(w)$.
\end{thm}

\subsubsection*{Global exponential stability for (\ref{is}),
infinite-horizon and invariant law bounds, and the CLT}

We state an important  global exponential stability result at $\tilde u$
for the non-linear dynamical system.
This is essential in the proof of the subsequent
infinite-horizon bounds for the marginals of the processes, which yield
 bounds on their long time limit, 
the invariant law.
We need  uniformity over the state space, and
results for the linearized equation~(\ref{linis}) 
are {\em not\/} enough. 

\begin{thm}
\label{glexst}
Let $\rho \le \theta < 1$
and
$u$ be the solution of (\ref{is}) starting at 
$u_0$ in $\mathcal{V} \cap L_2(g_\theta)$.
There exists  $\gamma_\theta>0$ and $C_\theta < \infty$
such that
$\Vert u_t -\ut \Vert_{L_2(g_\theta)} 
\le  \mathrm{e}^{- \gamma_\theta t} C_\theta \Vert u_0 -\ut \Vert_{L_2(g_\theta)}$.
\end{thm}

This does {\em not\/} hold in $L_2(\pi)$ for $L\ge2$, else 
Lemma~\ref{tight.in} below
would also hold in $L_2(\pi)$, which would contradict (\ref{finexp}).
Theorem 3.6 in Mitzenmacher~\cite{Mitzenmacher:96} states
a related result for some weighted $\ell_1$ norms obtained by potential 
function techniques. 

\begin{lem}
\label{tight.in}
Let  $\rho \le \theta < 1$ when $L\ge2$ or  $\rho < \theta < 1$ when $L=1$. Then
\[
\limsup_{N\to\infty} 
\EE\left(\left\Vert Z^N_0 \right\Vert_{L_2(g_\theta)}^2\right) < \infty
\Rightarrow
\limsup_{N\to\infty} \sup_{t\ge0}
\EE\left(\left\Vert Z^N_t \right\Vert_{L_2(g_\theta)}^2\right) < \infty
\]
and under the invariant laws
$
\limsup_{N\to\infty} \EE\left(\Vert Z^N_0 \Vert_{L_2(g_\theta)}^2\right) < \infty
$.
\end{lem}

Our main result is the functional CLT in equilibrium,
obtained with a compactness-uniqueness method
using tightness of the invariant laws (based on Lemma~\ref{tight.in}) 
and Theorems \ref{fclt} and \ref{ou.erg}.

\begin{thm}
\label{eqfclt}
Let the networks of size $N$ be in equilibrium.
For $L\ge2$ consider 
$L_2(g_{\rho})$
with its weak topology and
$\DD(\RR_+, L_2(g_{\rho}))$ with the corresponding Skorokhod
topology.
 Then $(Z^N)_{N\ge L}$
converges in law to the unique stationary OU~process 
solving the SDE~(\ref{sde}), in particular  $(Z^N_0)_{N\ge L}$ converges in law
to the invariant law for this process (see  Theorem~\ref{ou.erg}).
For $L=1$ the same result holds in  $L_2(g_{\theta})$ for $\rho<\theta<1$.
\end{thm}

\section{The proofs for converging initial conditions}
\setcounter{equation}{0}
\label{Sproofs} 

\subsection{Existence and uniqueness results}

\subsubsection*{Proof of Theorem \ref{isw} (refined existence result for (\ref{is}))}

We give the proof for $L_2(w)$, the proof for $\ell_1(w)$ being similar.
The assumption (\ref{compw}) and the identity 
$x^L-y^L = (x-y)(x^{L-1} + x^{L-2}y +\cdots + y^{L-1})$ yield
\begin{eqnarray*}
\left(u(k-1)^L - v(k-1)^L\right)^2 w(k)^{-1}
&\le& 
\left(u(k-1) - v(k-1)\right)^2 L^2 d w(k-1)^{-1}\,,
\\
\left(u(k)^L - v(k)^L\right)^2 w(k)^{-1}
&\le& 
\left(u(k) - v(k)\right)^2 L^2  w(k)^{-1}\,,
\\
\left(u(k+1) - v(k+1)\right)^2  w(k)^{-1}
&\le&
\left(u(k+1) - v(k+1)\right)^2  c^{-1}w(k+1)^{-1}\,,
\end{eqnarray*}
hence we have the Lipschitz bounds
$
\Vert F_+(u) - F_+(v) \Vert_{L_2(w)}^2 
\le 2\alpha^2 L^2 (d+1)
\Vert u - v  \Vert_{L_2(w)}^2
$ and
$
\Vert F_-(u) - F_-(v) \Vert_{L_2(w)}^2 
\le 2\beta^2  (c^{-1}+1)
\Vert u - v  \Vert_{L_2(w)}^2
$.
Existence and uniqueness follows by a classical
Cauchy-Lipschitz method.

\subsubsection*{The derivation of the Ornstein-Uhlenbeck process}

Let $(x)_k = x(x-1)\cdots(x-k+1)$ for $x\in \RR$ 
(the falling factorial of degree $k\in \NN$). Considering (\ref{F}),
let  the mappings $F^N_+$ and $F^N$ with values in $c^0_0$
be given for $v$ in $c_0$ by 
\begin{equation}
\label{FN}
F^N_+(v) (k) = \alpha\, {(Nv(k-1))_L - (Nv(k))_L\over (N)_L}\,,
\quad
k\ge1\,;
\qquad 
F^N(v) = F^N_+(v) - F^{\vphantom{N}}_-(v)\,.
\end{equation}
The process $R^N$ is  Markov on $\mathcal{V}^N$, and
when in state $r$ has jumps in its $k$-th coordinate, 
$k\ge1$, of size $1/N$  at rate $N F^N_+(r) (k)$
and size $-1/N$ at rate $N F_-(r) (k)$. 

\begin{lem}
\label{dyn}
Let 
$R^N_0$ be in  $\mathcal{V}^N$,
$u$ solve (\ref{is})
starting at $u_0$ in $\mathcal{V}$,
and $Z^N$ be given by (\ref{flu}). Then
\begin{equation}
\label{zeq}
Z^N_t  = Z^N_0
+ \int_0^t \sqrt{N}\left( F^N(R^N_s) - F(u_s) \right) ds  + M^N_t
\end{equation}
defines an independent family of square-integrable martingales
$M^N = (M^N(k))_{k\in\NN}$ 
independent of $Z^N_0$ with Doob-Meyer brackets given by
\begin{equation}
\label{bracket}
\left\langle M^N(k) \right\rangle_t 
= \int_0^t \left\{ F^N_+(R^N_s)(k) + F^{\vphantom{N}}_-(R^N_s)(k) \right\} ds\,.
\end{equation}
\end{lem}

\begin{proof}
This follows from a classical application of the 
Dynkin formula.
\end{proof}

The first lemma below shows that
it is indifferent to choose the
$L$ queues with or without replacement at this level of precision. The 
second one is a linearization formula.

\begin{lem}
\label{fact}
For $N\ge L\ge 1$ and $a$ in $\RR$ we have
\[
A^N(a) 
\defeg
{(Na)_L \over (N)_L} - a^L 
= \sum_{j=1}^{L-1} (a-1)^j a^{L-j}  
\sum_{1\le i_1<\cdots<i_j\le L-1}
{i_1\cdots i_j \over (N-i_1)\cdots (N-i_j)}
\]
and $A^N(a) =  N^{-1} O(a)$, uniformly for $0\le a \le 1$, and
$A^N(k/N) \le 0$ for $k=0$, $1$, \dots\,, $N$.
\end{lem}
\begin{proof}
We develop 
$
{(Na)_L \over (N)_L} = \prod_{i=0}^{L-1}  { Na-i \over N-i}
= \prod_{i=0}^{L-1} \left(a + (a-1){ i \over N-i}\right)
$
to obtain the identity for $A^N(a)$ which is clearly
$N^{-1} O(a)$, uniformly for $0\le a \le 1$.
For $a = k/N$, $\prod_{i=0}^{L-1}  { Na-i \over N-i}$ 
is composed of terms bounded by $a$
or contains a term equal to $0$ and cannot exceed $a^L$.
\end{proof}

\begin{lem}
\label{lin}
For $L\ge 1$ and $a$ and $h$ in $\RR$ we have
\[
B(a,h) 
\defeg
(a+h)^L - a^L - La^{L-1}h  
=\sum_{i=2}^{L} {L \choose i} a^{L-i} h^{i}
\]
with $B(a,h) =0$ for $L=1$ and $B(a,h) =h^2$
for $L=2$. For $L\ge2$ we have
$0\le B(a,h) \le h^{L} + \left(2^L -L -2\right) a h^2$
for $a$ and $a+h$ in $[0,1]$. 
\end{lem}
\begin{proof}
The identity is
Newton's binomial formula. 
A convexity argument yields $B(a,h)\ge 0$.
For $a$ and $a+h$ in $[0,1]$,
$
B(a,h) \le h^L + \sum_{i=2}^{L-1} {L \choose i} a h^2
= h^L + \left(2^L - L -2\right) a h^2\,.
$
\end{proof}

Let $v$ be in $\mathcal{V}$ and $x$ in $c_0^0$.
Considering (\ref{F}), (\ref{FN}) and Lemma~\ref{fact}, let
$G^N : \mathcal{V} \rightarrow  c^0_0$ be given  by
\begin{equation}
\label{defG}
G^N(v)(k)
=\alpha A^N(v(k-1)) -  \alpha  A^N(v(k))\,,
\qquad
k\ge1\,,
\end{equation}
and considering (\ref{F}), (\ref{defK}) and Lemma~\ref{lin}
let $H : \mathcal{V} \times c_0^0 \rightarrow  c_0^0$  be given by
\begin{equation}
\label{defH}
H(v,x)(k) = \alpha  B(v(k-1),x(k-1)) - \alpha B(v(k),x(k))\,,\qquad 
k\ge1\,,
\end{equation} 
so that for $v+x$ in $\mathcal{V}$ 
\begin{equation}
\label{diff}
F^N = F + G^N\,,
\qquad
F(v+x) - F(v)  
= \lin(v)  x + H(v,x)\,,
\end{equation}
and we derive the limit equation (\ref{ousde}) and (\ref{limdoob})
for the fluctuations from (\ref{zeq}) and (\ref{bracket}).

\subsubsection*{Proof of Theorem \ref{Kbdd} 
(existence and uniqueness for the OU~process)}

Considering (\ref{defK}), $v\le1$, convexity bounds, and (\ref{compw}),  we have 
\begin{eqnarray*}
\Vert \lin(v) x \Vert_{L_2(w)}^2
&\le&
2(\alpha L + \beta)
\sum_{k\ge1}
\left(
\alpha L x(k-1)^2 d w(k-1)^{-1} + 
(\alpha L + \beta) x(k)^2 w(k)^{-1}
\right.
\\
&&\kern 30mm\left. 
{} + \beta x(k+1)^2 c^{-1}  w(k+1)^{-1}
\right)
\\
&\le&
2(\alpha L + \beta)(\alpha L (d+1) + \beta (c^{-1} +1))
\Vert x \Vert_{L_2(w)}^2
\end{eqnarray*}
and (a) and (b) follow.
For $u_o$ in $\mathcal{V} \cap \ell_1(w)$
the martingale $M$ is square-integrable
in $L_2(w)$. If $\EE\Bigl(\Vert Z_0 \Vert_{L_2(w)}^2\Bigr) < \infty$ then 
the formula for $Z$ is well-defined, solves the SDE,
and the Gronwall Lemma yields 
$\EE\Bigl(\sup_{t\le T}\Vert Z_t \Vert_{L_2(w)}^2\Bigr) < \infty$. Else
for any $\varepsilon>0$
there is $r_\varepsilon<\infty$ such that 
$\PP\!\left(\Vert Z_0 \Vert_{L_2(w)} > r_\varepsilon\right) <\varepsilon$
and a localization procedure using pathwise uniqueness yields existence.

\subsection{The proof of the CLT}

\subsubsection*{Proof for Lemma \ref{l2bds} (finite-horizon bounds)}

Using (\ref{zeq}) and (\ref{diff})
\begin{equation}
\label{znavecg}
Z^N_t 
= Z^N_0 + M^N_t + \sqrt{N} \int_0^t G^N(R^N_s)\,ds
+\int_0^t 
\sqrt{N}\left(F (R^N_s) - F(u_s) \right)  ds
\end{equation}
where  Lemma~\ref{fact} yields 
$G^N(R^N_s)(k) = N^{-1} O\!\left( R^N_s(k-1)  +  R^N_s(k) \right)$ 
and considering (\ref{compw})
\begin{equation}
\label{grz}
\left\Vert G^N(R^N_s)\right\Vert_{L_2(w)}
=
N^{-1} O\!\left(\vphantom{R^N}\smash{\left\Vert R^N_s \right\Vert_{L_2(w)}}\right).
\end{equation}
We have
\begin{equation}
\label{decr1}
\left\Vert R^N_s \right\Vert_{L_2(w)}
\le
\left\Vert u_s \right\Vert_{L_2(w)} + 
N^{-1/2}  \left\Vert Z^N_s\right\Vert_{L_2(w)}
\end{equation}
and since $F_+$, $F_-$ and $F$ are Lipschitz (Theorem~\ref{isw})
the Gronwall Lemma 
yields that for some  $K_T<\infty$ we have
$
\left\Vert u_s \right\Vert_{L_2(w)}
\le K_T \left\Vert u_0 \right\Vert_{L_2(w)}
$ 
and 
\[
\sup_{0\le t\le T}\left\Vert Z^N_t \right\Vert_{L_2(w)}
\le
K_T\biggl(\left\Vert Z^N_0 \right\Vert_{L_2(w)} 
+ N^{-1/2} K_T \left\Vert u_0 \right\Vert_{L_2(w)}
+ \sup_{0\le t\le T}\left\Vert M^N_t \right\Vert_{L_2(w)}
\biggr).
\]
We conclude using the Doob inequality, (\ref{bracket}), (\ref{diff}),
(\ref{grz}), (\ref{decr1}), and
\begin{equation}
\label{brkbd}
\left\Vert F_+(R^N_s) + F_-(R^N_s) \right\Vert_{L_2(w)}
\le
K\left\Vert R^N_s \right\Vert_{L_2(w)}.
\end{equation}

\subsubsection*{Tightness for the process}

\begin{lem}
\label{tight}
Let $w$ satisfy (\ref{compw}), and consider 
$L_2(w)$
with its weak topology and
$\DD(\RR_+, L_2(w))$ with the corresponding Skorokhod
topology.
Let $u_0$ be in $\mathcal{V} \cap\ell_1(w)$ and
$R^N_0$  in $\mathcal{V}^N$, and
$Z^N$ be given by (\ref{flu}).
If $(Z^N_0)_{N\ge L}$ is tight
 then $(Z^N)_{N\ge L}$ is tight and its limit points are continuous.
\end{lem}

\begin{proof}
For $\varepsilon>0$ let $r_\varepsilon < \infty$ be such that
$\PP(Z^N_0 \in B(r_\varepsilon))> 1 - \varepsilon$ for $N\ge1$
(see the discussion prior to Theorem~\ref{fclt}).
Let $R^{N,\varepsilon}_0$ be equal to $R^N_0$ 
on $\{Z^N_0 \in B(r_\varepsilon)\}$ and such that
$Z^{N,\varepsilon}_0$ is uniformly
bounded in $L_2(w)$ on
$\{Z^N_0 \not\in B(r_\varepsilon)\}$.
Then $Z^{N,\varepsilon}_0$ is uniformly bounded in $L_2(w)$ and we may
 use a coupling argument
to construct $Z^{N,\varepsilon}$ and $Z^N$ coinciding
on $\{Z^N_0 \in B(r_\varepsilon)\}$.
Hence to prove tightness of $(Z^N)_{N\ge L}$ we may restrict our attention
to $(Z^N_0)_{N\ge L}$ uniformly bounded in $L_2(w)$, for which we may use 
Lemma \ref{l2bds}.

The compact subsets of $L_2(w)$ are Polish, a fact yielding tightness criteria. 
We deduce from Theorems~4.6 and 3.1 in 
Jakubowski~\cite{Jakubowski:86}, which considers
completely regular Hausdorff spaces (Tychonoff spaces)
of which $L_2(w)$ with its weak topology is an example,
that $(Z^N)_{N\ge L}$ is tight if
\begin{enumerate}
\item
For each $T\ge 0$ and $\varepsilon >0$ there is a bounded subset
$K_{T,\varepsilon}$ of $L_2(w)$ such that for
$N\ge L$ we have
$
\PP\!\left( Z^{N} \in \DD([0,T], K_{T,\varepsilon}) \right) > 1-\varepsilon
$.

\item
For each $d\ge1$, the 
$d$-dimensional processes $(Z^{N}(1), \ldots, Z^{N}(d))_{N \ge L}$
are tight.
\end{enumerate}

Lemma~\ref{l2bds} and the  Markov inequality yield condition 1.
We use (\ref{znavecg}) (see (\ref{zeq}) and (\ref{diff})),
 and (\ref{bracket}) and (\ref{diff}), and the  bounds (\ref{grz}), 
(\ref{decr1}) and (\ref{brkbd}). The bounds in Lemma~\ref{l2bds}
and the fact that $Z^{N}(k)$ has jumps of size $1/\sqrt{N} = o(N)$ classically imply
 that the above finite-dimensional processes are tight 
and have continuous limit points, 
see for instance  Ethier-Kurtz~\cite{Ethier:86} Theorem~4.1 p.~354
or Joffe-M{\'e}tivier~\cite{Joffe:86} Proposition~3.2.3 
and their proofs.
\end{proof}

\subsubsection*{Proof of Theorem \ref{fclt} (the functional CLT)}

Lemma \ref{tight} implies that from any subsequence of $Z^N$
we may extract a further subsequence which converges to
some $Z^\infty$ with continuous sample paths. Necessarily
$Z^\infty_0$ has same law as $Z_0$. In  (\ref{znavecg}) we have
considering (\ref{diff}) that
\begin{equation}
\label{fin}
\sqrt{N}\!\left(F (R^N_s)(k) - F(u_s)(k) \right)
= \lin(u_s) Z^N_s + \sqrt{N} H\!\left(u_s,\smash{N^{-1/2}}Z^N_s\right).
\end{equation}
We use the  bounds (\ref{grz}), (\ref{decr1}) and (\ref{brkbd}),
the uniform bounds in Lemma~\ref{l2bds},
and additionally (\ref{defH}) and Lemma~\ref{lin}. 
We deduce by a martingale characterization
that $Z^\infty$ has the law of the OU~process
unique solution for (\ref{ousde})
in $L_2(w)$
starting at $Z^\infty_0$, see Theorem~\ref{Kbdd}; 
the drift vector is given by the limit for (\ref{zeq}) and
(\ref{znavecg}) considering (\ref{fin}),
and the martingale bracket by the limit for (\ref{bracket}).
See for instance 
Ethier-Kurtz~\cite{Ethier:86} Theorem~4.1 p.~354 
or Joffe-M{\'e}tivier~\cite{Joffe:86} Theorem 3.3.1
and their proofs for details. Thus, this law is the unique
accumulation point for the relatively compact sequence of laws of $(Z^N)_{N\ge1}$,
which must then converge to it, proving Theorem \ref{fclt}.

\section{The properties of $\lint = \lin(\ut)$}
\setcounter{equation}{0}
\label{propOU}

\subsection{Proof of Theorem~\ref{Dbdd} (existence and uniqueness results)}

Considering (\ref{defKst}) and convexity bounds we have 
\begin{eqnarray*}
\Vert \lint z \Vert_{L_2(w)}^2
&=& \beta^2 \sum_{k\ge1}\left(
L \rho^{L^{k-1}} z(k-1)
-(L \rho^{L^{k}} + 1) z(k)
+  z(k+1)
\right)^2 w(k)^{-1}
\\
&\le&
\beta^2 (2L+2)\biggl(
L\sum_{k\ge1}
\rho^{2L^{k-1}} z(k-1)^2 w(k)^{-1}
+
L\sum_{k\ge1}
\rho^{2L^{k}} z(k)^2 w(k)^{-1}
\\
&&\kern30mm{}
+
\sum_{k\ge1} z(k)^2 w(k)^{-1}
+
\sum_{k\ge1} z(k+1)^2 w(k)^{-1}
\biggr)
\\
&\le&
\beta^2 (2L+2)\biggl(
L d \sum_{k\ge2} z(k-1)^2 w(k-1)^{-1}
+
(L \rho^{2L} +1) \sum_{k\ge1} z(k)^2 w(k)^{-1}
\\
&&\kern30mm{}
+
c^{-1}\sum_{k\ge1} z(k+1)^2 w(k+1)^{-1}
\biggr)
\\
&\le&
\beta^2 (2L+2) \left(L\rho^{2L} + Ld +c^{-1} +1\right) \Vert z \Vert_{L_2(w)}^2\,.
\end{eqnarray*}
The Gronwall Lemma yields uniqueness. For $k\ge1$ we have
\[
\left( L\rho^L \right)^{-1}  \pi(k+1) 
\le
\pi(k) =  \left( L\rho^{ L^{ \smash{k} } } \right)^{-1} \pi(k+1)  
\le \left( L^{-1} \rho^L \rho^{ -2L^{\smash{k}} } \right) \pi(k+1)\,,
\]
\[
\theta^{-1}\theta^{k+1} 
\le  \theta^k 
\le  \left( \theta^{-1} \rho^L  \rho^{ -2L^{\smash{k}} } \right) \theta^{k+1} \,.
\]
When $B$ is an Hilbertian Brownian motion,
the formula for $Z$ yields a well-defined solution.

\subsection{A related birth and death process, and the spectral decomposition}
\label{spectdec}

Considering (\ref{matrixK}), $\mathcal{A} = \lint^*$ is  
the infinitesimal generator   
of the sub-Markovian birth and death process on  the irreducible class $(1,2,\ldots)$
with birth rates $\lambda_k =\beta L \rho^{L^k}$ 
and death rates $\mu_k=\beta$
for $k\ge1$ (killed at rate $\mu_1=\beta$ at state $1$).
The process is well-defined since the rates are bounded.

Karlin and McGregor~\cite{Karlin:57a,Karlin:57b} give a spectral 
decomposition for such  processes, used by
Callaert and Keilson~\cite{Callaert:73a,Callaert:73b} 
and van Doorn~\cite{Doorn:85}
to study exponential ergodicity
properties.
The state space in these works is $(0,1,2,\dots)$,
possibly extended by an absorbing barrier 
or graveyard state at
$-1$. We consider  $(1,2,\ldots)$ and adapt their notations to this simple shift.

The potential coefficients 
(\cite{Karlin:57a} eq.~(2.2), \cite{Doorn:85} eq.~(2.10)) are
given  by 
\[
\pi(k) 
= {\lambda_1 \cdots \lambda_{k-1} \over \mu_2 \cdots \mu_k}
=  L \rho^{L^1}\cdots  L \rho^{L^{k-1}} =   L^{k-1} \rho^{(L^k-L)/(L-1)},
\qquad k\ge1\,,
\]
and solve the detailed balance
equations $\mu_{k+1} \pi(k+1) = \lambda_{k} \pi(k)$ with
$\pi(1)=1$, see (\ref{pot}). 

The equation ${\mathcal A}Q(x) = - x Q(x)$
for an eigenvector $Q(x) = (Q_n(x))_{n\ge1}$ of eigenvalue $-x$ 
yields $\lambda_1 Q_2(x) = (\lambda_1 + \mu_1 -x) Q_1(x)$ and
$\lambda_n Q_{n+1}(x) = (\lambda_n + \mu_n -x) Q_n(x) - \mu_n Q_{n-1}(x)$ for $n\ge2$.
With the natural convention $Q_0=0$
and normalizing choice $Q_1=1$, we obtain
inductively  $Q_n$ as the  polynomial of degree $n-1$ satisfying
the recurrence relation
\[
- x  Q_n(x)
=  \beta  Q_{n-1}(x)  -  \left(\beta L \rho^{L^n}   + \beta  \right) Q_n(x) 
+ \beta  L \rho^{L^n}  Q_{n+1}(x)\,,
\qquad  n\ge1\,,
\]
corresponding to 
\cite{Karlin:57a}~eq.~(2.1)
and \cite{Doorn:85}~eq.~(2.15). Such 
a sequence of polynomials is orthogonal with respect to a probability measure
$\psi$ on $\RR_+$ and, for $i,j\ge1$ with $i\neq j$,
$
\int_0^\infty Q_i(x)^2\, \psi(dx) = \pi(i)^{-1}
$ 
and
$
\int_0^\infty Q_i(x) Q_j(x)\, \psi(dx) = 0
$
or in matrix notation
$\int_0^\infty Q(x) Q(x)^*\, \psi(dx) = \mathrm{diag}(\pi^{-1})$.

Let $P_t = (p_t(i,j))_{i,j\ge1}$ denote the sub-stochastic 
transition matrix for $\mathcal{A}$. The adjoint matrix
$P_t^*$ 
is the fundamental solution for the forward equation 
$\dot z_t = \mathcal{A}^* z_t = \lint z_t$ given in (\ref{linis}).
The representation formula of 
Karlin and McGregor~\cite{Karlin:57a,Karlin:57b}, see
(1.2) and (2.18) in \cite{Doorn:85},
yields
\begin{equation}
\label{repr}
\mathrm{e}^{\lint t} = P_t^* = (p^*_t(i,j))_{i,j\ge1}\,,
\qquad
p^*_t(i,j) = p_t(j,i) = 
\pi(i) \int_0^\infty 
\mathrm{e}^{-xt} Q_i(x) Q_j(x)\, \psi(dx)\,,
\end{equation}
or in  matrix notation
$\mathrm{e}^{\lint t} 
= \mathrm{diag}(\pi) \int_0^\infty  \mathrm{e}^{-xt} Q(x) Q(x)^* \,\psi(dx)$.

The probability measure $\psi$ is called the spectral measure, its support $S$
is called the spectrum, and we set  $\gamma=\min S$. The 
OU~process in Theorem~\ref{Dbdd}~(b) and its invariant law
and its covariance matrix in Theorems~\ref{ou.erg} and \ref{eqfclt} 
can be written
\begin{eqnarray}
\label{explou}
Z_t &=&
\mathrm{diag}(\pi)
\int_S  \mathrm{e}^{-xt}   Q(x)^*  \left(Z_0
+\int_0^t\mathrm{e}^{xs} \,dB_s 
\right) Q(x) \,\psi(dx)\,,
\\
\label{limlaw}
\int_0^\infty \mathrm{e}^{\lint t}\,dB_t
&=&
\mathrm{diag}(\pi) \int_S 
\left(Q(x)^*\int_{0}^\infty  \mathrm{e}^{-xt} \,dB_t \right) Q(x) \,\psi(dx)\,,
\\
\label{cov}
\int_0^\infty \mathrm{e}^{\lint t} \mathrm{diag}(\tilde v) \mathrm{e}^{\lint^* t}\,dt
&=&
\mathrm{diag}(\pi)
\int_{S^2} 
{Q(x)^* \mathrm{diag}(\tilde v) Q(y) \over x+y}\,
Q(x) Q(y)^*  \,\psi(dx)\psi(dy)\,
\mathrm{diag}(\pi).\qquad
\end{eqnarray}

\subsection{The spectral gap, exponential stability, and ergodicity}
\label{prf.spe.gap}

\subsubsection*{Proof of Theorem~\ref{spe.gap} 
(spectral gap and exponential stability in the self-adjoint case)}

The potential coefficients $(\pi(k))_{k\ge1}$ solve the detailed balance equations 
for $\mathcal{A}$ and hence 
$\lint = \mathcal{A}^*$ is self-adjoint in 
$L_2(\pi)$.
For the spectral gap,
we follow Van Doorn~\cite{Doorn:85}, Section~2.3.
The orthogonality properties imply that 
 $Q_n$ has $n-1$ distinct zeros $0 < x_{n,1}<\ldots<x_{n,n-1}$ such that
$x_{n+1,i} <x_{n,i} < x_{n+1,i+1}$
for $1\le i \le n-1$. Hence $\xi_i = \lim_{n\to\infty }x_{n,i} \ge 0$
exists, $\xi_i\le \xi_{i+1}$, and $\sigma = \lim_{i\to\infty}\xi_i$ exists
in $[0,\infty]$. 
Theorem 5.1 in \cite{Doorn:85} establishes that
$\gamma>0$ if and only if $\sigma>0$ and
Theorem~5.3~(i) in \cite{Doorn:85} that
$\sigma = \left(\sqrt {\lim_k \lambda_k} 
- \sqrt {\lim_k \mu_k}  \right)^2 = \beta>0$.
(Theorem 3.3 in \cite{Doorn:85} states that $\gamma = \xi_1 \le \sigma$, but
estimating $\xi_1$ is impractical.)

For the exponential stability, we have
$\Vert z_t \Vert_{L_2(\pi)}^2 = \left(\mathrm{e}^{\lint t} z_0, 
\mathrm{e}^{\lint t} z_0 \right)_{L_2(\pi)}$ and
the fact that
$\mathrm{e}^{\lint t}$ is self-adjoint in $L_2(\pi)$ and 
the spectral representation (\ref{repr}) yield
\begin{eqnarray*}
\left(\mathrm{e}^{\lint t} z_0, 
\mathrm{e}^{\lint t} z_0 
\right)_{L_2(\pi)}
&=&
\left( z_0, 
\mathrm{e}^{ 2\lint t} z_0 
\right)_{L_2(\pi)}
=
\int_S \mathrm{e}^{- 2xt} z_0^* Q(x) Q(x)^* z_0\,\psi(dx) 
\\
&\le& \mathrm{e}^{- 2\gamma t}  
\int_S z_0^* Q(x) Q(x)^* z_0\,\psi(dx) 
=
\mathrm{e}^{- 2\gamma t} \left( z_0,z_0 \right)_{L_2(\pi)}.
\end{eqnarray*}

\subsubsection*{Proof of Theorem~\ref{linis.exp} (exponential 
stability, non self-adjoint case)}

It is similar to and simpler than the proof for Theorem~\ref{glexst}
to which Section~\ref{expstab} is devoted,
and we postpone the proof until the end of that section.

\subsubsection*{Proof of Theorem~\ref{ou.erg} 
(ergodicity for the OU~process)}

We use the uniqueness result and
explicit formula in Theorem~\ref{Dbdd},
and Theorem~\ref{spe.gap} or \ref{linis.exp}.

\section{Exponential stability for the nonlinear system}
\setcounter{equation}{0}
\label{expstab}

\subsection{Some comparison results}

Considering  (\ref{diff}) with $\lint = \lin(\ut)$  and $F(\ut)=0$,
if $u$ solves (\ref{is}) in $\mathcal{V}$
then $y = u-\ut$ solves the 
recentered equation  given by 
$\dot y_t (k)  = F(\tilde u + y) = \lint y_t(k) + H(\ut, y_t)(k)$ or
\begin{eqnarray}
\label{center}
\kern-5mm
\dot y_t (k) 
&=&
\beta L \rho^{L^{k-1}}y_t(k-1) + \alpha  B(\ut(k-1),y_t(k-1))  
\nonumber\\
&&\kern7mm{}
-\left(\beta L \rho^{L^k}y_t(k) + \alpha B(\ut(k),y_t(k)) + \beta y_t(k) \right) 
+  \beta y_t(k+1)\,,
\qquad k\ge1\,.
\end{eqnarray}
If $u_0$ is in $\mathcal{V}\cap \ell_1$ then $u$ is in $\mathcal{V}\cap \ell_1$
and hence $y$ is in $\ell_1^0$ and for $k\ge1$ we have
\begin{equation}
\label{sum.cen}
\dot y_t(k) + \dot y_t(k+1) + \cdots\,
=  \beta L \rho^{L^{k-1}} y_t(k-1) + \alpha B(\ut(k-1), y_t(k-1))  - \beta y_t(k)\,.
\end{equation}
If $y$ solves (\ref{center}) starting at 
$y_0$ such that $y_0+ \ut$ is in $\mathcal{V}$, then $u = y+ \ut$ 
solves (\ref{is}) in $\mathcal{V}$
starting at $u_0=y_0+ \ut$. Then $-\ut \le y \le 1-\ut$ and $-1 < y  < 1$.
For $y_0+ \ut$ in $\mathcal{V}\cap\ell_1$ we have  $y$ in $\ell_1^0$. 

\begin{lem}
\label{order}
Let $u$ and $v$ be two solutions for (\ref{is}) 
in $\mathcal{V}$ such that $u_0 \le v_0$. Then $u_t \le v_t$ for $t\ge0$.
Let $y_0+ \ut$ be in $\mathcal{V}$ and $y$ solve (\ref{center}).
If $y_0\ge0$ then $y_t\ge0$ and if $y_0\le0$ then $y_t\le0$ for $t\ge0$.
\end{lem}

\begin{proof}
Lemma 6  in \cite{Vved:96} yields the result for  (\ref{is}) (the proof written 
for $L=2$ is valid for $L\ge1$).
The result for (\ref{center}) follows by considering $u=y+\ut$ and $\ut$ 
which solve (\ref{is}).
\end{proof}

We compare solutions of the nonlinear equation (\ref{center}) 
and of certain linear equations.

\begin{lem}
\label{kolm}
Let $\hat\mathcal{A}$ be the generator of the 
sub-Markovian birth and death process
with birth rate $\hat\lambda_k \ge 0$ 
and death rate $\beta$ at $k\ge1$. Let $\sup_k \hat\lambda_k<\infty$.
The linear operator $x\mapsto \hat\mathcal{A}^*x$ given by
\[
\hat\mathcal{A}^*x(k) = 
\hat\lambda_{k-1} x(k-1)
-(\hat\lambda_{k} + \beta ) x(k)
+  \beta x(k+1)\,,
\qquad
k\ge1\,,
\]
is bounded in $\ell_1^0$.
There exists a unique $z=(z_t)_{t\ge0}$ given by 
$z_t = \mathrm{e}^{\hat\mathcal{A}^* t} z_0$
solving
the forward  Kolmogorov equation $\dot z=\hat\mathcal{A}^*z$ in $\ell_1^0$. It is 
such that
if $z_0\ge0$ then $z_t\ge0$ and if
$z_0\le0$ then $z_t\le0$, and 
$
\dot z_t(k) + \dot z_t(k+1) + \cdots\,
=   \hat\lambda_{k-1} z_t(k-1)  - \beta z_t(k)
$  for $k\ge1$.
\end{lem}

\begin{proof}
The  operator norm in $\ell_1^0$ of $\hat\mathcal{A}^*$ 
is bounded by $2(\sup_k \hat\lambda_k + \beta)$, hence existence and uniqueness.
Uniqueness and linearity imply that if $z_0=0$ then $z_t=0$
and else if
$z_0\ge0$ then 
$z_t \Vert z_0\Vert_1^{-1}$ is the instantaneous
law of the process starting at  $z_0\Vert z_0\Vert_1^{-1}$
and hence $z_t\ge0$. If $z_0\le0$ then $- z$ solves the equation
starting at $- z_0 \ge 0$ and hence $- z_t \ge 0$. The last result is obtained by 
summation.
\end{proof}

\begin{lem}
\label{cen.lin}
Let $L\ge2$ and $y=(y_t)_{t\ge0}$ solve (\ref{center}) with  $y_0+ \ut$ in 
$\mathcal{V}\cap\ell_1$.
Under the assumptions  of Lemma~\ref{kolm},
let $z=(z_t)_{t\ge0}$  solve 
$\dot z  = \hat\mathcal{A}^*z $  in $\ell_1^0$. Let
$h=(h_t)_{t\ge0}$ be given in $\ell_1^0$ by
\[
h(k) = z(k) + z(k+1) + \cdots\,
- (y(k) + y(k+1) + \cdots\,)\,,
\qquad
k\ge1\,.
\]

\noindent
(a) Let $\hat\lambda_k \ge \beta L \rho^{L^k} 
+  \alpha\!\left(1 + \left(2^L-L-2\right) \ut(k)\right)$  for $k\ge1$,
$y_0\ge0$, and $h_0\ge0$. Then $h_t\ge0$ for $t\ge0$.

\noindent
(b) Let $\hat\lambda_k \ge \beta L \rho^{L^k}$  for $k\ge1$,
$y_0\le0$, and $h_0\le0$.   Then $h_t\le0$ for $t\ge0$.
\end{lem}

\begin{proof}
We prove (a).
For $\varepsilon >0$ let $\hat\mathcal{A}^*_\varepsilon$ correspond to
$\hat\lambda_k^\varepsilon = \hat\lambda_k + \varepsilon$.
The operator norm in $\ell_1^0$ of
$\hat\mathcal{A}^*_\varepsilon -\hat\mathcal{A}^*$ is bounded by
$2 \varepsilon$, hence
$\lim_{\varepsilon\to0} \mathrm{e}^{\hat\mathcal{A}_\varepsilon^* t} z_0 =z_t$ 
in $\ell_1^0$ and we may 
assume that 
$\hat\lambda_k > \beta L \rho^{L^k}
+  \alpha\!\left(1 + \left(2^L-L-2\right) \ut(k)\right)$ for $k\ge1$.
Since
$z_t = \mathrm{e}^{\hat\mathcal{A}^* t} z_0$ depends continuously on $z_0$ 
in $\ell_1^0$
we may assume $h_0 > 0$. 
Let $\tau = \inf\{t\ge0 : \{k\ge1 : h_{t}(k)=0\} \neq\emptyset \}$  be
the first time when $h(k)=0$ for some $k\ge1$. 
We have $\tau>0$.

The result (a) holds if  $\tau =\infty$. 
If $\tau \not=\infty$, 
Lemma \ref{kolm} and (\ref{sum.cen}) yield
\begin{eqnarray*}
\dot h_\tau(k) &=&
\hat\lambda_{k-1}y_\tau(k-1)  - \beta L \rho^{L^{k-1}} y_\tau(k-1)
- \alpha B(\ut(k-1), y_\tau(k-1))
\nonumber\\
&&\kern10mm{}
+\hat\lambda_{k-1} (z_\tau(k-1) - y_\tau(k-1) )
-\beta  (z_\tau(k) -  y_\tau(k))\,.
\end{eqnarray*}
Lemma~\ref{order} yields $y \ge 0$ and 
Lemma~\ref{lin}  and $y\le1$ yield
\begin{eqnarray*}
B(\ut(k-1), y(k-1))&\le&
y(k-1)^L + \left(2^L -L -2\right)\ut(k-1) y(k-1)^2
\\
&\le& \left(1 + \left(2^L -L -2\right)\ut(k-1)\right) y(k-1)\,,
\end{eqnarray*}
hence
$
\hat\lambda_{k-1}y(k-1)  - \beta L \rho^{L^{k-1}} y(k-1)
- \alpha B(\ut(k-1), y(k-1)) \ge0
$
with equality only when $y(k-1) = 0$.
For $k$ in $\mathcal{Z} = \{k\ge1 : h_\tau(k)=0 \} \neq \emptyset$ we have
\[
z_\tau(k-1) - y_\tau(k-1) =  h_\tau(k-1)\ge 0\,,
\quad
z_\tau(k) - y_\tau(k)  = - h_\tau(k+1) \le 0 \,,
\]
hence $\dot h_\tau(k) \ge 0$
with equality if only if 
$k-1$ is in $\mathcal{Z}\cup\{0\}$ and $k+1$ is in $\mathcal{Z}$.
Moreover $h_t(k) > 0$ for $t<\tau$
and $h_\tau(k) = 0$ imply $\dot h_\tau(k) \le 0$. Hence
$\dot h_\tau(k) = 0$, and the above signs and equality cases  yield that
$z_\tau(k-1)=y_\tau(k-1)=0$ and
$k-1$ is in $\mathcal{Z}\cup\{0\}$ and $k+1$ is in $\mathcal{Z}$.
By induction $z_\tau(i) = y_\tau(i)=0$ for $i\ge1$ which implies 
$z_{t} = y_{t}=0$ for $t\ge \tau$, and the proof of (a) is complete.

The proof for (b) is similar and involves obvious changes of sign.
The assumption $\hat\lambda_k > \beta L \rho^{L^k}$ suffices to conclude since
$B(\ut(k-1), y(k-1))\ge0$ (Lemma~\ref{lin}) and
the non-linearity ``pushes'' in the right direction.
\end{proof}

\begin{lem}
\label{equiv}
For any $0<\theta<1$ there exists $K_\theta<\infty$ such that
for $x$ in $L_2(g_\theta)\subset \ell^0_1$ 
\[
\left\Vert (x(k)+x(k+1)+\cdots)_{k\ge1} \right \Vert_{L_2(g_\theta)}
\le
 K_\theta \Vert x \Vert_{L_2(g_\theta)}\,.
\]
\end{lem}

\begin{proof}
Using a classical convexity inequality
\begin{eqnarray*}
&&\kern-2mm
\sum_{k\ge1}  (x(k)+ x(k+1) + \cdots\, )^2 \theta^{-k} 
\\
&&\kern2mm {} \le
\sum_{k\ge1} n\!\left(
x(k)^2 + x(k+1)^2 + \cdots + x(k+n-2)^2 
+ (x(k+n-1)+ x(k+n) + \cdots\, )^2
\right) \theta^{-k} 
\\
&&\kern2mm {} \le 
n \!\left(1+\theta + \cdots+ \theta^{n-2} \right)\sum_{k\ge1}  x(k)^2 \theta^{-k}
+ n\,  \theta^{n-1}\sum_{k\ge1}  (x(k)+ x(k+1) + \cdots\, )^2 \theta^{-k}
\end{eqnarray*} 
and we take $n$ large enough that $n\theta^{n-1}<1$ and
$
K_\theta^2 =  n \!\left(1+\theta + \cdots+ \theta^{n-2} \right)
(1 - n \theta^{n-1})^{-1}
$.
\end{proof}

\subsection{Proofs of the exponential stability results}

\subsubsection*{Proof of Theorem~\ref{glexst} for $L\ge2$}

If  $u_0$ is in $\mathcal{V} \cap L_2(g_\theta)$, then so are
$u^-_0 = \min\{u_0,\ut\}$ and $u^+_0 = \max\{u_0,\ut\}$
and hence the 
corresponding solutions $u^-$ and $u^+$ for (\ref{is}), see
Theorem~\ref{isw}.
Lemma~\ref{order} yields that $u^-_t  \le u_t \le u^+_t$ and 
$u^-_t  \le \ut \le u^+_t$ 
for $t\ge0$. Then
\[
y= u-\ut\,,
\qquad
y^+ = u^+ - \ut \ge 0\,,
\qquad
y^- =  u^- - \ut \le 0\,,
\]
solve the recentered equation (\ref{center}), and termwise
\begin{equation}
\label{trmin}
|y_0| = \max\{y^+_0 , - y^-_0\}\,,
\qquad
|y_t| \le \max\{y^+_t , - y^-_t\} \,,
\quad t\ge0\,.
\end{equation}

We consider the birth and death process with generator $\hat\mathcal{A}$
defined in Lemma~\ref{kolm}
with 
\[
\hat \lambda_k 
= \max\left\{ 
\beta L \rho^{L^{k}} +   \alpha\!\left(1 + \left(2^L-L-2\right) \ut(k)\right) , 
\beta \theta 
\right\}\,,
\qquad
k\ge1\,,
\]
which satisfies the assumptions of Lemma~\ref{cen.lin} (a) and (b).
We reproduce the spectral study in Section~\ref{spectdec} 
and the proof of Theorem~\ref{spe.gap} in  Section~\ref{prf.spe.gap}
for $\hat\mathcal{A}$, 
corresponding objects being denoted with a hat.
For $\rho \le \theta <1$ we have $\alpha \le \beta \theta$ and hence
$\hat \lambda_k$ is equivalent to $\beta \theta$ for large $k$,
Theorems 5.1 and 5.3~(i) in \cite{Doorn:85} yield that
$0<\hat\gamma  \le \hat\sigma 
= \left(\sqrt{\beta\theta} - \sqrt{\beta}\right)^2 
= \beta \left(1 - \sqrt{\vphantom{\beta}\theta} \right)^2$, and 
if $z$ solves $\dot z = \hat\mathcal{A}^* z$ then
$\Vert z_t \Vert_{L_2(\hat\pi)}
\le \mathrm{e}^{-\hat\gamma t}  \Vert z_0 \Vert_{L_2(\hat\pi)}$
for $t\ge0$. Moreover 
\[
\theta^{k-1} 
\le
\hat\pi(k) = \theta^{k-1} 
\prod_{i=1}^{k-1}
\max\left\{ 
\theta^{-1}L \rho^{L^{k}} 
+ \theta^{-1}\rho\!\left(1 + \left(2^L-L-2\right) \ut(k)\right),
1 \right\} 
\]
and the product converges, hence
$\hat\pi(k) = O(\theta^{k})$ and $\theta^{k} = O(\hat\pi(k))$ and
Lemma~\ref{eqiv.l2w} yields that there exists
$c>0$ and $d>0$ such that 
$c^{-1} \Vert \cdot \Vert_{L_2(\hat\pi)} \le 
\Vert \cdot \Vert_{L_2(g_\theta)} \le d \Vert \cdot \Vert_{L_2(\hat\pi)}$.
Hence for $t\ge0$
\[
\Vert z_t \Vert_{L_2(g_\theta)} 
\le d \Vert z_t \Vert_{L_2(\hat\pi)}
\le \mathrm{e}^{-\hat\gamma t}  d  \Vert z_0 \Vert_{L_2(\hat\pi)}
\le  \mathrm{e}^{-\hat\gamma t} cd \Vert z_0 \Vert_{L_2(g_\theta)}\,.
\]

Hence if $z^+$ solves $z^+= \hat\mathcal{A}^* z^+$ starting at $z_0^+ =y_0^+ \ge 0$
then Lemmas~\ref{cen.lin} (a) and \ref{equiv} yield
\begin{eqnarray*}
\Vert y^+_t \Vert_{L_2(g_\theta)} 
&\le &
\Vert (y^+_t(k) + y^+_t(k+1) +\cdots\, )_{k\ge1} \Vert_{L_2(g_\theta)}
\\
&\le &
\Vert (z^+_t(k) +  z^+_t(k+1) +\cdots\, )_{k\ge1} \Vert_{L_2(g_\theta)}
\\
&\le &
K_\theta \Vert z^+_t \Vert_{L_2(g_\theta)} 
\le \mathrm{e}^{-\hat\gamma t} cd K_\theta \Vert y^+_0 \Vert_{L_2(g_\theta)}
\end{eqnarray*}
and similarly 
if $z^-$ solves $z^-= \hat\mathcal{A}^* z^-$ starting at $z_0^- =y_0^- \le0$
then Lemmas~\ref{cen.lin} (b) and \ref{equiv} yield
$
\Vert y^-_t \Vert_{L_2(g_\theta)} \le  
\mathrm{e}^{-\hat\gamma t} cd K_\theta \Vert y^-_0 \Vert_{L_2(g_\theta)}
$.
We set $\gamma_\theta = \hat\gamma$ and $C_\theta = cd K_\theta$.
Considering (\ref{trmin}),
\[
\Vert y_t \Vert_{L_2(g_\theta)}^2 
\le
\Vert y^+_t \Vert_{L_2(g_\theta)}^2 + \Vert y^-_t \Vert_{L_2(g_\theta)}^2 
\le \mathrm{e}^{- 2\gamma_\theta t} C_\theta^2
\left(
\Vert y^+_0 \Vert_{L_2(g_\theta)}^2 + \Vert y^-_0 \Vert_{L_2(g_\theta)}^2 
\right)
\]
and 
we complete the proof by remarking that
for $k\ge1$, either $y^+_0(k) = y_0(k)$ and  $y^-_0(k) = 0$ or
$y^-_0(k) = y_0(k)$ and $y^+_0(k) = 0$, and hence
$
\Vert y^+_0 \Vert_{L_2(g_\theta)}^2 + \Vert y^-_0 \Vert_{L_2(g_\theta)}^2 
 = 
\Vert y_0 \Vert_{L_2(g_\theta)}^2$.

\subsubsection*{Proof of Theorem~\ref{linis.exp} and 
of Theorem~\ref{glexst} for  $L=1$}

The linearization (\ref{linis}) of Equation (\ref{is}) 
is obtained by replacing $B$ and $H$ in Equation (\ref{center})  by $0$
and
coincides with Equation (\ref{center}) for $L=1$.
 Likewise,
the equation for (\ref{linis}) corresponding to (\ref{sum.cen})
is obtained by omitting the terms $\alpha B(\ut(k-1), y_t(k-1))$.
We obtain a result for Equation  (\ref{linis}) 
corresponding to Lemma~\ref{cen.lin} (a) and (b)
under the sole assumption
$\hat \lambda_k \ge \beta L \rho^{L^{k}}$ for $k\ge1$.
The proof proceeds as for Theorem~\ref{glexst} for $L\ge2$ 
with the difference that
$\hat \lambda_k = \max\left\{\beta L \rho^{L^{\smash{k}}}, \beta \theta \right \}$.
We have $\hat \lambda_k$  equal to $\beta \theta$ for large $k$
for $0<\theta< 1$ when  $L\ge2$ and for $\rho \le \theta <1$
when $L=1$.

\section{Tightness estimates and the functional CLT in equilibrium}
\label{cfp}
\setcounter{equation}{0}

\subsection{Proof of Lemma~\ref{tight.in} (infinite horizon and
invariant law bounds)}

Let  $U_h (v)$ be the solution of (\ref{is}) at time 
$h\ge0$ with initial value $v$ in $\mathcal{V}$. For $t_0\ge0$ let
$
Z^N_{t_0,h} = \sqrt{N}\left(R^N_{t_0+h}-  U_h(R^N_{t_0}) \right).
$
Then 
$Z^N_{t_0+h} = Z^N_{t_0,h} + \sqrt{N}\left( U_h(R^N_{t_0}) - \ut \right) $
and Theorem~\ref{glexst} yields 
\begin{equation}
\label{decineg}
\left\Vert Z^N_{t_0+h}  \right\Vert_{L_2(g_\theta)}
\le 
\left\Vert
Z^N_{t_0,h}
\right\Vert_{L_2(g_\theta)}
+
\mathrm{e}^{- \gamma_\theta h} C_\theta 
\left\Vert
Z^N_{t_0}
\right\Vert_{L_2(g_\theta)}.
\end{equation}

The conditional law of
$(Z^N_{t_0,h})_{h\ge0}$ 
given $R^N_{t_0} =r$
is the law of $Z^N$
started with $R^N_0 = u_0 =r$,
in particular with $Z^N_0=Z^N_{t_0,0}=0$.
We reason as in 
(\ref{znavecg})--(\ref{brkbd}) 
except that the bound  (\ref{decr1}) 
becomes 
\[
\left\Vert R^N_{t_0+s} \right\Vert_{L_2(g_\theta)}
\le
\left\Vert \tilde u \right\Vert_{L_2(g_\theta)} + 
N^{-1/2}  \left\Vert Z^N_{t_0 +s} \right\Vert_{L_2(g_\theta)}
\]
and we use (\ref{decineg})  and  obtain that for some $K_T<\infty$
\[
\sup_{0\le h \le T} \left\Vert Z^N_{t_0,h} \right\Vert_{L_2(g_\theta)} 
\le K_T  \biggl(
 N^{-1/2} \left\Vert \ut \right\Vert_{L_2(g_\theta)}
+
N^{-1}  C_\theta 
\left\Vert
Z^N_{t_0}
\right\Vert_{L_2(g_\theta)}
+
\sup_{0\le h \le T} \left\Vert M^N_{t_0+h} - M^N_{t_0} \right\Vert_{L_2(g_\theta)}
\biggr)
\]
which combined with (\ref{decineg}) yields
that for some $L_T<\infty$ we have for $0\le h \le T$
\begin{equation}
\label{bonineg}
\EE\left(
\left\Vert Z^N_{t_0+h}  \right\Vert_{L_2(g_\theta)}^2
\right)
\le 
L_T
+
2 ( K_T N^{-1} + \mathrm{e}^{- \gamma_\theta h})^2 C_\theta^2 
\,\EE\left(
\left\Vert
Z^N_{t_0}
\right\Vert_{L_2(g_\theta)}^2
\right).
\end{equation}

We fix $T$ large enough for 
$8\mathrm{e}^{- 2\gamma_\theta T}C_\theta^2  \le \varepsilon < 1$.
Uniformly for $N \ge K_T \mathrm{e}^{\gamma_\theta T}$, 
for $m\in\NN$
\[
\EE\left(
\left\Vert Z^N_{\smash{(m+1)}T}  \right\Vert_{L_2(g_\theta)}^2
\right)
\le 
L_T
+
\varepsilon
\,\EE\left(
\left\Vert
Z^N_{mT}
\right\Vert_{L_2(g_\theta)}^2
\right)
\]
and  by induction
\[
\EE\left(
\left\Vert Z^N_{mT}  \right\Vert_{L_2(g_\theta)}^2
\right)
\le 
L_T
\sum_{j=1}^m
\varepsilon^{j-1}
+
\varepsilon^m
\,\EE\left(
\left\Vert
Z^N_{0}
\right\Vert_{L_2(g_\theta)}^2
\right)
\le 
{L_T \over 1 - \varepsilon} + 
\EE\left(\left\Vert
Z^N_{0}
\right\Vert_{L_2(g_\theta)}^2
\right),
\]
and (\ref{bonineg}) also yields
\[
\sup_{0 \le h \le T}
\EE\left(
\left\Vert Z^N_{mT+h}  \right\Vert_{L_2(g_\theta)}^2
\right)
\le 
L_T
+
8 C_\theta^2 
\,\EE\left(
\left\Vert
Z^N_{mT}
\right\Vert_{L_2(g_\theta)}^2
\right),
\]
hence the infinite horizon bound
\[
\sup_{t\ge0}\EE\left(
\left\Vert Z^N_{t}  \right\Vert_{L_2(g_\theta)}^2
\right) 
\le
L_T
+
8 C_\theta^2 
\left(
{L_T \over 1 - \varepsilon} + 
\EE\left(\left\Vert
Z^N_{0}
\right\Vert_{L_2(g_\theta)}^2
\right)
\right).
\]

Ergodicity
and the Fatou Lemma yield that for $Z_\infty^N$ 
distributed according to the invariant law
\[
\EE\left(\left\Vert Z^N_\infty \right\Vert_{L_2(g_\theta)}^2\right)
\le
\liminf_{t\ge0} 
\EE\left(\left\Vert Z^N_t \right\Vert_{L_2(g_\theta)}^2\right)
\le
\sup_{t\ge0} 
\EE\left(\left\Vert Z^N_t \right\Vert_{L_2(g_\theta)}^2\right)
\]
and the invariant law bound
follows if we show that we can choose $R^N_0$ in $\mathcal{V}^N$ such
that
\begin{equation}
\label{initgood}
\limsup_{N\to\infty} 
\EE\left(\left\Vert Z^N_0 \right\Vert_{L_2(g_\theta)}^2\right) < \infty\,.
\end{equation}
For this
we consider $L\ge2$, the case $L=1$ being similar, and $R_0^N$
given for $k\ge0$ by
$R_0^N(k) = iN^{-1}$ with $1\le i \le N$ such that
$- 2^{-1} N^{-1} < \ut (k) - iN^{-1} \le  2^{-1} N^{-1}$.
For $x\ge 0$ and $0< y\le 1$
\begin{eqnarray*}
y = \rho^{(L^x-1)/(L-1)}
&\Leftrightarrow&
x = \log\left(1 + (L-1)\log y / \log \rho \right) / \log L
\\
&\Leftrightarrow&
\theta^{-x} = \left(1 + (L-1) \log y /\log \rho \right)^{-\log\theta / \log L}
\end{eqnarray*}
hence for
$
z(N) =\inf\!\left\{k\ge 1 : R_0^N(k)=0 \right\}$
we have
$
z(N)
= \inf\!\left\{k\ge 1 : \ut (k) \le 2^{-1} N^{-1} \right\}
= \inf\left\{k\in\NN 
: k \ge \log\left(1 + (L-1)\log \left( 2^{-1} N^{-1} \right) / \log \rho \right) 
/ \log L
\right\}
$.
Then
\[
\left\Vert Z^N_0 \right\Vert_{L_2(g_\theta)}^2
=
N \sum_{k=1}^{z(N)-1} \left(R^N_0(k) - \ut (k)\right)^2 \theta^{-k}
+
N \sum_{k\ge z(N)} \ut (k)^2 \theta^{-k}
\]
with
\[
N \sum_{k=1}^{z(N)-1} \left(R^N_0(k) - \ut (k)\right)^2 \theta^{-k}
\le 2^{-2}N^{-1}\,{\theta^{-z(N)}-\theta^{-1}  \over \theta^{-1} -1 }
= O \left( N^{-1} (\log N)^{\smash{-\log \theta / \log L}}  \right)
\]
and for large enough $N$ (and hence $z(N)$)
\begin{eqnarray*}
N \sum_{k\ge z(N)} \ut (k)^2 \theta^{-k}
&=&
N \ut (z(N))^2 \,\sum_{j\ge0 } \rho^{2L^{z(N)}(L^j-1)/ (L-1)}\theta^{-(j+z(N))}
\\
&\le&
2^{-2}N^{-1}\,\sum_{j\ge0 } \rho^{L^{z(N)}(L^j-1)/ (L-1)} = o(N^{-1})\,,
\end{eqnarray*}
hence (\ref{initgood}) holds and the proof is complete.

\subsection{The functional CLT: Proof of Theorem~\ref{eqfclt}}
\label{pclt}

Lemma~\ref{tight.in} and the Markov inequality imply that
in equilibrium
$(Z_0^N)_{N\ge L}$ is  tight for the weak
topology of $L_2(g_\rho)$, for which all bounded sets are relatively compact.
Consider a subsequence.
We can extract a further 
subsequence along which $(Z_0^N)_{N\ge L}$ converges 
in law to some square-integrable
$Z_0^\infty$
in $L_2(g_\rho)$, and
Theorem~\ref{fclt} yields that along the further subsequence $(Z^{N})_{N\ge L}$
converges in law to the OU~process $Z^\infty$
unique solution for (\ref{sde})
in $L_2(g_{\rho})$
starting at $Z^\infty_0$. 

The limit in law of a sequence of stationary processes is stationary
(Ethier-Kurtz~\cite{Ethier:86} p.~131, Lemma~7.7 and Theorem~7.8).
Hence the law of $Z^\infty$ is determined as the unique law of the stationary 
OU~process given by (\ref{sde}), see Theorem~\ref{ou.erg}.
From every subsequence we can extract a further
subsequence converging in law to $Z^\infty$, hence
$\lim_{N\to\infty} Z^N = Z^\infty$ in law.

\vskip5mm
\noindent
{\em Acknowledgment\/}.
The author would like to thank referees and editors for their valuable 
suggestions, in particular of writing this combined version replacing
two precedent preprints which separated the equilibrium and 
non-equilibrium studies.



\end{document}